\documentclass{article}
\usepackage{amssymb}
\usepackage{mathrsfs}
\usepackage{cite}

\usepackage{stmaryrd}
\usepackage{amsfonts}
\usepackage{amsmath}
\usepackage{bm}
\usepackage{indentfirst}
\usepackage{array}
\usepackage{arydshln}
\usepackage{graphicx}
\usepackage{listings}
\usepackage{epstopdf}
\numberwithin{equation}{section}
\usepackage{geometry}
\usepackage{amsthm,color}
\newtheorem{theorem}{Theorem}

\newtheorem{lemma}{Lemma}

\newtheorem{rem}{Remark}
\newtheorem{proposition}{Proposition}
\title{\bf{Stability of viscous shock profile for convective porous-media flow with degenerate viscosity$^\ast$}
\setcounter{footnote}{-1}
\author{Yechi Liu$^\dag$
\noindent\footnote{
$\ast$\quad This research is supported by the Natural Science Foundation of Hunan Province, China (Grant No. 2023JJ40659).\\
\indent $\,\,\,\,\,\dag$\quad E-mail: lyc9009@sina.cn.}\\
\small College of Science, National University of Defence Technology, Changsha 410003, P.R.China}}
\date{}
\begin{document}\large
\maketitle

\textbf{Abstract}. In this paper, we are concerned with the large time behavior of viscous shock wave
for the convective porous-media equation with degenerate viscosity. We get the regularity of the solution for general initial data and prove the shock wave is nonlinearly stable providing the initial perturbation is small. Moreover, the $L^\infty$ decay rate is obtained, which generalized the famous result \cite{osh82}. Note that the traditional energy method and continuity argument can not be directly used in this paper since the degeneration of viscosity. One need to fully utilize the sign of perturbation and it derivatives, decompose the integral domain to ensure that in each domain the sign is invariant. Then the stability and the decay rate are obtained by energy method and an area inequality.

%{It was proved  in \cite{osh82} that there exists viscous shock wave for the convective porous-media equation with degenerate viscosity. In this paper, we get the regularity of the solution and prove such viscous shock wave is nonlinearly stable and also obtain the converges rate in $L^\infty$-norm as time goes to infinity. The proof is given by the elementary energy method.\\}

\textbf{Keywords}. Porous-media flow, asymptotic behavior, degenerate viscosity, viscous shock wave, decay rate.

\section{Introduction and main results}
We are concerned with the quasi-linear parabolic equation
\begin{equation}\label{001}
u_t+f(u)_x=A(u)_{xx},\qquad t>0,\quad x\in\mathbb R,
\end{equation}
where $f,A\in C^2$. In addition, $a(u)=:A^\prime(u)>0$ for any $u\neq0$ and $a(0)=0$. When $A(u)=u^m\,(m>1)$ and $f(u)=-u^n\,(n\in\mathbb N_+)$, \eqref{001} becomes the convective porous-media equation, and the existence, regularity and finite propagating speed of solutions were proved in Gilding-Peletier \cite{gil76} and Gilding \cite{gil77}. We refer to \cite{lad67,gil96,dib93} for general $f(u)$ and $A(u)$.

Similar to the Burgers equation, the equation \eqref{001} also admits viscous shock waves. In fact, it was proved in \cite{osh82} that there exists viscous shock wave of \eqref{001}. Furthermore, the authors proved the $L^1$ stability of viscous shock waves under the assumption that the initial values stay between far field end states. This condition was subsequently relaxed by Freist\"{u}hler-Serre \cite{fre98} for the linear diffusion case, i.e., $A(u)=u$, and by Feireisl-Lauren\c{c}ot \cite{fei99} for porous-media type. See also a nice survey \cite{ser04}.

In the case of $A(u)=u$, there are remarkable works considering the decay rate of shock wave, see \cite{KM} and the reference thereafter. Especially, Nishihara-Zhao \cite{nis02} further obtained the convergence rate toward the viscous shock waves in $L^\infty$-norm under some restriction conditions on the initial data. Kang-Vasseur \cite{kan17} showed similar results in $L^2$-norm without such restrictions on initial data. Huang-Xu \cite{hua22} obtained the decay rate without assuming that the initial perturbation belongs to some weighted Sobolev space. We also refer to \cite{k,I,HX2,LWX} and the references therein for the decay rates of the rarefaction wave and multi-dimensional case.

While for the case we studied in this paper, there are much more less researches about the decay rate of the shock wave. The main difficulty comes from the degeneracy of viscosity. because of which , the equation will be perfect nonlinear when considering the anti-derivative, so that we can not use the classical energy method directly. In this paper, we firstly prove the existence of the solution and obtain the boundedness of the derivative. Then, we separate the whole integral domain into intervals with a standard of the sign about the perturbation and its derivative, and estimate the energy function case by case. At last, we use the area inequality to obtain the decay rate of the perturbation in $L^2-$norm, and hence in $L^\infty-$norm.

\vskip 0.2in

Now we give the main theorem. In this paper, we study the Cauchy problem of \eqref{001} for porous media type $A(u)=u^m,1<m<2$, that is,
\begin{equation}\label{0}
\left\{\begin{aligned}
&u_t+f(u)_x=(u^m)_{xx},\\
&u(0,x)=u_0(x).
\end{aligned}\right.
\end{equation}
Here, $f$ represents the flux function with $f(0)=0$ and $f^{\prime\prime}\geqslant C_f>0$. It is noted that the equation \eqref{0} is degenerate parabolic when $u=0$ since $m>1$. Due to the physical meaning of the problem, we assume that $u\geqslant0$. Then, we have the following global existence for general initial data.
\begin{theorem}\label{thm1}
Assume $1<m<2$, $0\leqslant u_0(x)\in L^\infty(\mathbb R)\cap C(\mathbb R)$, and $u_0^m,u_0^{m-1}$ are Lipschitz continuous. The Cauchy problem \eqref{0} and \eqref{upm} admits a global in time solution $u(t,x)$ satisfies\\
\begin{equation*}
\left\{\begin{aligned}
&u(t,x)\in L^\infty\big([0,\infty)\times\mathbb R\big)\cap C^{\frac{1}{2},1}\big((0,\infty)\times\mathbb R\big),\\
&u_x(t,\cdot),\big(u^m\big)_x(t,\cdot)\in L^\infty\big((0,\infty)\times\mathbb R\big),\\
&u_x(t,\cdot),\big(u^m\big)_x(t,\cdot)\in C(\mathbb R)
\end{aligned}\right.
\end{equation*}
for any $t>0$.
\end{theorem}

\vskip 0.2in

If additionally,
\begin{equation}\label{upm}
\lim_{x\rightarrow-\infty}u_0(x)=u_->0,\qquad \lim_{x\rightarrow+\infty}u_0(x)=u_+:=0,
\end{equation}
there exists a corresponding viscous shock wave of \eqref{0}
\begin{equation}\label{U}
u(t,x)=U(\xi),\qquad \xi=:x-\gamma t
\end{equation}
satisfying $\lim_{\xi\rightarrow\pm\infty}U(\xi)=u_\pm$, where $\gamma$ is a constant given by the Rankine-Hugoniot condition
\begin{equation}\label{rhc}
f(u_+)-f(u_-)=\gamma(u_+-u_-).
\end{equation}
Hence, it is difficult to use the standard energy estimate to study the asymptotic behavior of the solution, or the stability of the viscous shock wave.

\vskip 0.2in

Given the viscous shock wave $U(x)$ mentioned in \eqref{U}, if $u_0(x)-U(x)\in L^1(\mathbb R)$, there exists a space shift $x_0$ satisfying
\begin{equation}\label{phi0}
\int_\mathbb R\big(u_0(x)-U(x+x_0)\big)\textrm{d}x=0.
\end{equation}
Without loss of generality, we take $x_0=0$ in what follows. Now we can state the time-decay rate as follows.
\begin{theorem}\label{thm2}
Let $1<m\leqslant\frac{4}{3}$ and let $u(t,x)$ be the solution given in Theorem \ref{thm1}. In addition, assume $u_0-U\in L^1(\mathbb R)$, $\Phi_0=:\int_{-\infty}^\xi\big(u_0(\eta)-U(\eta)\big)\textrm{d}\eta\in L^2(\mathbb R)$. Then there exists a small constant $\varepsilon_0>0$ such that, when
\begin{equation}\label{005}
\|\Phi_0\|_{H^1(\mathbb R)}\leqslant\varepsilon_0,
\end{equation}
the solution $u(t,x)$ in (1) satisfies
\begin{equation}\label{jie1}
\|u(t,\cdot)-U(\cdot-\gamma t)\|_{L^2(\mathbb R)}\leqslant C_\delta(1+t)^{-\frac{1}{4(11m+7)}+\delta},
\end{equation}
where $\delta$ is any small positive constant and $C_\delta>0$ is a constant depending on $\delta$.
\end{theorem}
\begin{rem}
Since $u_x\in L^\infty\big((0,\infty)\times\mathbb R\big)$ from Theorem \ref{thm1} and $U^\prime$ is bounded from Remark \ref{Urem} given in Section 2, we can conclude from \eqref{jie1} that
\begin{equation*}
\|u(t,\cdot)-U(\cdot-\gamma t)\|_{L^\infty(\mathbb R)}\leqslant C_\delta(1+t)^{-\frac{1}{6(11m+7)}+\delta}
\end{equation*}
by using the interpolation inequality.
\end{rem}

\vskip 0.2in

The rest of this paper is organized as follows. In Section 2, we will give some properties about the viscous shock wave $U$ and derive the perturbation equation. Then in Section 3, the proof on the existence and regularity of the solution (i.e. Theorem 1) is given. At last, the time decay rate (i.e. Theorem 2) will be obtained in Section 4.

\vskip 0.2in

\noindent\textbf{Notations}. For function spaces, $L^p=L^p(\mathbb R)$ and $H^k=H^k(\mathbb R)$ denote the usual Lebesgue space and $k-$th order Sobolev space on the whole space $\mathbb R$ with norms $\|\cdot\|_p$ and $\|\cdot\|_{H^k}$, respectively, which means
\begin{equation*}
\|v\|_p=:\left(\int_\mathbb R|v(x)|^p\,\textrm{d}x\right)^\frac{1}{p},\quad \|v\|_{H^k}=:\left(\sum_{l=0}^k\|\partial_x^lv\|_2^2\right)^\frac{1}{2},
\end{equation*}
where $\partial_x^lv=\frac{\partial^lv}{\partial x^l}$. We also denote $\|\cdot\|=\|\cdot\|_2$ for simplicity. For the functions with space shift, we denote
\begin{equation*}
g_{(y)}(x)=:g(x-y)
\end{equation*}
for any $x,y\in\mathbb R$. We also use $c$ and $C$ to represent uncertain positive constants suitably small and large respectively.

\section{Preliminaries}
Firstly, we recall the viscous shock wave of \eqref{0} constructed in \cite{osh82}. Denote the viscous shock wave by
\begin{equation}\label{99vs}
U(\xi)=:U(x-\gamma t),
\end{equation}
which satisfies $\lim_{\xi\rightarrow\pm\infty}U(\xi)=u_\pm$, $u_->u_+=0$. Then we have the following Lemma.
\begin{lemma}[\!\cite{osh82}\,]\label{hhhh}
Let $U(\xi)$ be the viscous shock wave given by \eqref{99vs}, then $U\in C^1(\mathbb R)$ and $U^\prime\leqslant0$. Furthermore, if $m=1$, $U(\xi)>0$ for all $\xi\in\mathbb R$; if $m>1$, there exist some $x_R$ such that $U(\xi)=0$ for all $\xi\geqslant x_R$ and $U(\xi)>0$ for all $\xi<x_R$.
\end{lemma}
\begin{rem}\label{Urem}
In fact, $U$ satisfies
\begin{equation*}
mU^{m-1}U^\prime=f(U)-f(u_-)-\gamma(U-u_-),
\end{equation*}
which implies
\begin{equation*}
f^\prime(0)-\gamma\leqslant U^\prime\leqslant0.
\end{equation*}
\end{rem}

\vskip 0.2in

We recall some properties of solutions to parabolic equation \eqref{0}. Denote the solution semigroup of \eqref{0} as $T(t)$, it holds that
\begin{equation*}
T(t)u_0(x)=u(t,x),\qquad x\in\mathbb R,t\geqslant0,
\end{equation*}
then we have
\begin{lemma}[\!\cite{cra80,osh82}\,]\label{lema} $T(t)$ has the following properties
\begin{equation*}
\begin{aligned}
&{\rm (1)}\quad T(t)\:\;{\rm commutes\,\, with\,\, translation\!:}\,\,T(t)u_{(y)}=\big(T(t)u\big)_{(y)};\\
&{\rm (2)}\quad T(t)\:\;{\rm is\,\, monotone\!:}\,\,u_0(x)\leqslant v_0(x)\Rightarrow\big(T(t)u_0\big)(x)\leqslant\big(T(t)v_0\big)(x), a.e.;\\
&{\rm (3)}\quad T(t)\:\;{\rm preserves}\,\,L^1\!: u_0-v_0\in L^1(\mathbb R)\Rightarrow T(t)u_0-T(t)v_0\in L^1(\mathbb R);\\
&{\rm (4)}\quad T(t)\:\;{\rm is\,\, conservative\!:}\,\,u_0-v_0\in L^1(\mathbb R)\Rightarrow\int_\mathbb R\big(T(t)u_0-T(t)v_0\big){\rm d}x=\int_\mathbb R(u_0-v_0){\rm d}x;\\
&{\rm (5)}\quad T(t)\:\;{\rm is\,\, contractive\,\, in}\,\,L^1\!: \|T(t)u_0-T(t)v_0\|_1\leqslant\|u_0-v_0\|_1.
\end{aligned}
\end{equation*}
\end{lemma}
With the help of Lemma \ref{lema}, we can conclude from \eqref{phi0} that, for any $t\geqslant0$,
\begin{equation}\label{alpha}
\int_\mathbb R\big(u(t,x)-U(x-\gamma t)\big)\textrm{d}x=0.
\end{equation}
Note that we have supposed the space shift $x_0$ to be $0$. Define the perturbation $\phi(t,\xi)=u(t,\xi+\gamma t)-U(\xi)$. Since $\phi(t,\cdot)\in L^1(\mathbb R)$ for any $t>0$ and is uniformly continuous with respect to $x$, we can conclude that $\lim_{\xi\rightarrow\pm\infty}\phi(t,\xi)=0$. Thus, $\phi$ satisfies
\begin{equation}\label{phi}
\left\{\begin{aligned}
&\phi_t-\gamma\phi_\xi+\big(f(U+\phi)-f(U)\big)_\xi=\big((U+\phi)^m-U^m\big)_{\xi\xi},\\
&\phi(0,\xi)=u_0(\xi)-U(\xi)=:\phi_0(\xi),\\
&\lim_{\xi\rightarrow\pm\infty}\phi(t,\xi)=0.
\end{aligned}\right.
\end{equation}
Owing to \eqref{alpha}, we can define
\begin{equation}\label{phixx}
\Phi(t,\xi)=\int_{-\infty}^\xi u(t,\eta+\gamma t)-U(\eta)\textrm{d}\eta
\end{equation}
which satisfies
\begin{equation}\label{Phi}
\left\{\begin{aligned}
&\Phi_t-\gamma\Phi_\xi+f(U+\Phi_\xi)-f(U)=\big((U+\Phi_\xi)^m-U^m\big)_\xi,\\
&\Phi(0,\xi)=\Phi_0(\xi),\\
&\lim_{\xi\rightarrow\pm\infty}\Phi(t,\xi)=0.
\end{aligned}\right.
\end{equation}

\vskip 0.2in

To deal with the viscosity, the following inequalities are needed.
\begin{proposition}[Lemma 4.4 on page 13 in \cite{dib93}]
Suppose $a,b\in\mathbb R$ and $\mu\geqslant1$, it holds that
\begin{equation}\label{ab}
|a-b|^{\mu+1}\leqslant C_\mu(|a|^{\mu-1}a-|b|^{\mu-1}b)(a-b)
\end{equation}
for some constant $C_\mu>0$ depending only on $\mu$.
\end{proposition}
\begin{proposition}
Suppose $a,b\geqslant0$ and $0<\mu\leqslant1$, it holds that
\begin{equation*}
\left|a^\mu-b^\mu\right|\leqslant C_\mu|a-b|^\mu
\end{equation*}
for some constant $C_\mu>0$ depending only on $\mu$.
\end{proposition}
\begin{proposition}\label{phii}
For any $p\geqslant2,1<m\leqslant\frac{4}{3}$ and $w(x)\in H^1(\mathbb R)$ satisfying $w_x\in L^{m+1}(\mathbb R)$, it holds
\begin{equation}\label{103a}
\int_\mathbb R|w|^{p-1}w_x^2\textrm{d}x\leqslant C\|w\|_\frac{2-m}{m-1}^{2-m}\int_\mathbb R|w|^{p-2}|w_x|^{m+1}\textrm{d}x.
\end{equation}
\end{proposition}
\begin{proof}
With the H\"{o}lder's inequality, it holds
\begin{equation}\label{104}
\int_\mathbb R|w|^{p-1}w_x^2\textrm{d}x\leqslant\|w\|_{\kappa_1}^{\kappa_1\frac{m-1}{m+1}}\left(\int_\mathbb R|w|^{p-2}|w_x|^{m+1}\textrm{d}x\right)^{\frac{2}{m+1}},
\end{equation}
where $\kappa_1$ is a positive constant satisfying
\begin{equation}\label{kap1}
\kappa_1=p-1+\frac{2}{m-1}>\frac{2-m}{m-1}\geqslant2.
\end{equation}
On the other hand, the interpolation inequality implies that
\begin{equation}\label{106}
\begin{aligned}
\|w\|_{2\kappa_1}^{\frac{m+p-1}{m+1}}&=\left\||w|^{\frac{m+p-1}{m+1}}\right\|_{2\kappa_1\frac{m+1}{m+p-1}}\\
&\leqslant C\left\||w|^{\frac{p-2}{m+1}}w_x\right\|_{m+1}^{\kappa_2}
    \left\||w|^{\frac{m+p-1}{m+1}}\right\|_{\frac{2-m}{m-1}\frac{m+1}{m+p-1}}^{1-\kappa_2}\\
&=C\left(\int_\mathbb R|w|^{p-2}|w_x|^{m+1}\textrm{d}x\right)^{\frac{\kappa_2}{m+1}}\|w\|_\frac{2-m}{m-1}^{\frac{m+p-1}{m+1}(1-\kappa_2)},
\end{aligned}
\end{equation}
where $\kappa_2\in(0,1)$ satisfies
\begin{equation}\label{kap2}
\kappa_2\left(\frac{1}{m+1}-1\right)+(1-\kappa_2)\frac{m-1}{2-m}\frac{m+p-1}{m+1}=\frac{m+p-1}{2\kappa_1(m+1)}.
\end{equation}
Furthermore, noting that $\kappa_1>2$, and using the H\"{o}lder's inequality and \eqref{106}, we have
\begin{equation}\label{107}
\begin{aligned}
\|w\|_{\kappa_1}^{\kappa_1}&\leqslant\|w\|_\frac{2-m}{m-1}^{\kappa_3}\|w\|_{2\kappa_1}^{\kappa_1-\kappa_3}\\
&\leqslant C\|w\|_\frac{2-m}{m-1}^{\kappa_3+(\kappa_1-\kappa_3)(1-\kappa_2)}\left(\int_\mathbb R|w|^{p-2}|w_x|^{m+1}\textrm{d}x\right)^{\frac{\kappa_1-\kappa_3}{m+p-1}\kappa_2},
\end{aligned}
\end{equation}
where $\kappa_3\in(0,\kappa_1)$ is a constant satisfying
\begin{equation}\label{kap3}
\kappa_3\frac{m-1}{2-m}+\frac{\kappa_1-\kappa_3}{2\kappa_1}=1.
\end{equation}
Comparing \eqref{104} and \eqref{107}, we complete the proof.
\end{proof}
\begin{proposition}\label{phiii}
For any $p>2,1<m<2$ and $w(x)\in H^1(\mathbb R)\cap W^{1,\infty}(\mathbb R)$, it holds that
\begin{equation}\label{402a}
\left(\|w\|_p^p\right)^\nu\leqslant C\int_\mathbb R|w|^{p-2}|w_x|^{m+1}\textrm{d}x,
\end{equation}
where
\begin{equation*}
\nu=1+\frac{3m+1}{p-2}.
\end{equation*}
\end{proposition}
\begin{proof}
With the help of the interpolation inequality, we have
\begin{equation}\label{ew1}
\begin{aligned}
\|w\|_\infty^{\frac{m+p-1}{m+1}}&=\left\||w|^{\frac{m+p-1}{m+1}}\right\|_\infty\\
&\leqslant C\left\||w|^{\frac{p-2}{m+1}}w_x\right\|_{m+1}^\nu
    \left\||w|^{\frac{m+p-1}{m+1}}\right\|_{\frac{m+1}{m+p-1}p}^{1-\nu}\\
&=C\left(\int_\mathbb R|w|^{p-2}|w_x|^{m+1}\textrm{d}x\right)^{\frac{\nu}{m+1}}\|w\|_p^{\frac{m+p-1}{m+1}(1-\nu)},
\end{aligned}
\end{equation}
where $\nu=\frac{m+p-1}{mp+m+p-1}$. Since $p>2$, H\"{o}lder's inequality and \eqref{ew1} imply that
\begin{equation*}
\begin{aligned}
\|w\|_p^p&\leqslant\|w\|^2\|w\|_\infty^{p-2}\leqslant C\|w\|_\infty^{p-2}\\
  &\leqslant C\left(\int_\mathbb R|w|^{p-2}|w_x|^{m+1}\textrm{d}x\right)^{\frac{\nu}{m+p-1}(p-2)}\|w\|_p^{(1-\nu)(p-2)},
\end{aligned}
\end{equation*}
which completes the proof.
\end{proof}
In addition, the following lemma will be used to get the time-decay rate of $\phi$.
\begin{lemma}[Lemma 2.3 in \cite{hua22}]
Assume $f(t)\in C^1[0,\infty)\cap L^1[0,\infty)$ to be any non-negative function satisfying
\begin{equation*}
\frac{{\rm d}\!f}{{\rm d}t}\leqslant (1+t)^{-\alpha},\qquad 0<\alpha\leqslant2.
\end{equation*}
Then, it holds that
\begin{equation*}
f(t)\leqslant C(1+t)^{-\frac{\alpha}{2}}.
\end{equation*}
\end{lemma}

\section{Existence and Regularity}

In this section, we will prove Theorem \ref{thm1}. The proof is separated into several lemmas. Firstly, we have
\begin{lemma}\label{48}
Assume $u_0$ is continuous in $\mathbb R$, $0\leqslant u_0\leqslant M$ and $u_0^m$ is Lipschitz continuous, then the Cauchy problem \eqref{0} admits a bounded continuous weak solution $u(t,x)$ on $(0,T]\times\mathbb R$ for any constant $T>0$ and satisfies that $(u^m)_x$ is bounded and $0\leqslant u\leqslant M$. Furthermore, $u(t,x)$ is a classical solution on $\{(t,x)\mid u(t,x)>0\}$.
\end{lemma}
\begin{proof}
Denote $v_0=u_0^m$. From the assumption, we can construct a sequence of smooth functions $\{v_{0,n}(x)\}$ which uniformly converges to $v_0(x)$ and satisfies
\begin{equation*}
\left|\frac{\textrm{d}}{\textrm{d}x}v_{0,n}\right|\leqslant K,\quad \frac{1}{n}\leqslant v_{0,n}\leqslant M^m,\quad n=1,2,\cdots,
\end{equation*}
where $K>2M^m$ is a constant. Choose a sequence of truncation $\{w_n(x)\}$ satisfying
\begin{equation}\label{pr1}
\left\{\begin{aligned}
&w_n(x)=v_{0,n}(x),\qquad |x|\leqslant n-2,\\
&w_n(x)=M^m,\qquad\quad\, |x|\geqslant n-1,\\
&\frac{1}{n}\leqslant w_n\leqslant M^m,\quad \left|\frac{\textrm{d}}{\textrm{d}x}w_n\right|\leqslant K, \quad n=1,2,\cdots.
\end{aligned}\right.
\end{equation}
Let $v=u^m$ and define $\alpha(v)=v^{\frac{1}{m}}$, then \eqref{0} becomes
\begin{equation}\label{pr2}
\alpha^\prime(v)v_t=-f\big(\alpha(v)\big)_x+v_{xx}.
\end{equation}
Consider the initial boundary value problem of \eqref{pr2} with the following initial boundary data
\begin{equation}\label{pr3}
v(0,x)=w_n(x),\quad v(t,\pm n)=M^m.
\end{equation}
From the theory of classical parabolic equation, the problem \eqref{pr2} and \eqref{pr3} has a classical solution $v_n(t,x)$ satisfying
\begin{equation}\label{zz0}
\frac{1}{n}\leqslant\inf w_n(x)\leqslant v_n(t,x)\leqslant M^m.
\end{equation}

\vskip 0.2in

\noindent (i).We will then prove $\frac{\partial}{\partial x}v_n$ is uniformly bounded on $\Omega_n=[0,T]\times[-n,n]$ with respect to $n$.

Let $P_n=\frac{\partial}{\partial x}v_n$, then from \eqref{pr2}, $P_n$ satisfies
\begin{equation*}
\alpha^\prime(v_n)\frac{\partial}{\partial t}P_n=\frac{\partial^2}{\partial x^2}P_n-\left(\frac{\big(\alpha^\prime(v_n)\big)_x}{\alpha^\prime(v_n)}
+f^\prime\big(\alpha(v_n)\big)\alpha^\prime(v_n)\right)\frac{\partial}{\partial x}P_n-f^{\prime\prime}\big(\alpha(v_n)\big)\alpha^\prime(v_n)^2P_n^2.
\end{equation*}
Using the maximum principle (Theorem 2.9 on page 23 in \cite{lad67}), it holds
\begin{equation*}
\max_{\Omega_n}|P_n|\leqslant\max_{\Gamma_n}|P_n|,
\end{equation*}
where $\Gamma_n$ is the parabolic boundary of $\Omega_n$. On $t=0$, from \eqref{pr1} we have
\begin{equation*}
\left|\frac{\partial}{\partial x}v_n\right|=\left|\frac{\textrm{d}}{\textrm{d}x}w_n\right|\leqslant K.
\end{equation*}
On $x=n$, using the maximum principle, it holds $v_n(t,n)=\max_{\Omega_n}v_n$. Thus,
\begin{equation*}
\left.\frac{\partial}{\partial x}v_n\right|_{x=n}\geqslant0.
\end{equation*}
Let $z_n=v_n-M^m(x-n+1)$. Then $z_n$ satisfies
\begin{equation*}
\alpha^\prime(v_n)z_{nt}=-f\big(\alpha(v_n)\big)_x+z_{nxx}
\end{equation*}
on $Q_n=[0,T]\times[n-1,n]$. It is easy to see that $z_n(0,x)\geqslant0$ for $x\in[n-1,n]$ and $z_n(t,n)=0,z_n(t,n-1)>0$ for $t\in[0,T]$. Thus, $z_n(n,t)=\min_{Q_n}z_n$, which implies
\begin{equation*}
\left.\frac{\partial}{\partial x}z_n\right|_{x=n}\leqslant0.
\end{equation*}
Hence,
\begin{equation*}
\left|\frac{\partial}{\partial x}v_n\right|_{x=n}\leqslant M^m.
\end{equation*}
Similarly, it holds $\left|\frac{\partial}{\partial x}v_n\right|_{x=-n}\leqslant M^m$. Using the maximum principle, we obtain
\begin{equation}\label{pr4}
\left|\frac{\partial}{\partial x}v_n\right|\leqslant\max\{K,M^m\}=K.
\end{equation}

\vskip 0.1in

\noindent (ii). We will prove that $v_n$ is uniformly H\"{o}lder continuous with respect to $n$ and the index is $\left\{\frac{1}{2},1\right\}$. Choose $n$ sufficiently large. For any $t\in[0,T]$ and $x_1,x_2\in[-n,n]$, we have
\begin{equation}\label{pr5}
|v_n(t,x_1)-v_n(t,x_2)|\leqslant K|x_1-x_2|.
\end{equation}
Let $u_n=\alpha(v_n)$. From \eqref{pr2}, it holds
\begin{equation}\label{pr6}
\frac{\partial}{\partial t}u_n=-\frac{\partial}{\partial x}f(u_n)+\frac{\partial^2}{\partial x^2}v_n.
\end{equation}
For any $s,t\in[0,T]$, denote $\Delta t=t-s$ and since $n$ is large, we can ensure that $x,x+|\Delta t|^{\frac{1}{2}}$ are both in $[-n,n]$. Integrating \eqref{pr6} over $[s,t]\times[x,x+|\Delta t|^{\frac{1}{2}}]$ implies
\begin{equation*}
\begin{aligned}
&\left|\int_x^{x+|\Delta t|^{\frac{1}{2}}}\big(u_n(t,y)-u_n(s,y))\textrm{d}y\right|\\
&=\left|\int_s^t\left(f\big(u_n(\tau,x)\big)-f\big(u_n(\tau,x+|\Delta t|^{\frac{1}{2}})\big)+\frac{\partial}{\partial x}v_n(\tau,x+|\Delta t|^{\frac{1}{2}})-\frac{\partial}{\partial x}v_n(\tau,x)\right)\textrm{d}\tau\right|\\
&\leqslant C|\Delta t|.
\end{aligned}
\end{equation*}
Using the mean value theorem for integral, there exists a $x^\ast\in[x,x+|\Delta t|^{\frac{1}{2}}]$ such that
\begin{equation*}
|u_n(t,x^\ast)-u_n(s,x^\ast)|\leqslant C|\Delta t|^{\frac{1}{2}}.
\end{equation*}
Thus,
\begin{equation*}
|v_n(t,x^\ast)-v_n(s,x^\ast)|\leqslant mM^{m-1}|u_n(t,x^\ast)-u_n(s,x^\ast)|\leqslant C|\Delta t|^{\frac{1}{2}},
\end{equation*}
where $C$ is independent of $n$. This inequality, together with \eqref{pr5} implies that, for any $(t,x),(s,y)\in[0,T]\times[-n,n]$, when $n$ is sufficiently large, it holds
\begin{equation*}
\begin{aligned}
&|v_n(t,x)-v_n(s,y)|\\
&\leqslant |v_n(t,x)-v_n(t,x^\ast)|+|v_n(t,x^\ast)-v_n(s,x^\ast)|+|v_n(s,x^\ast)-v_n(s,y)|\\
&\leqslant C(|t-s|^{\frac{1}{2}}+|x-y|),
\end{aligned}
\end{equation*}
where $C$ is independent of $n$.

\vskip 0.1in

\noindent (iii).These conclusions above imply that the sequence $\{u_n\}\big(u_n=\alpha(v_n)\big)$ is uniformly bounded and equicontinuous. Then, from the Arzela-Ascoli theorem, there exists a subsequence, still denote as it self, converges uniformly on any compact subset of $[0,T]\times\mathbb R$. Since $f\in C^2$ and $\{u_n\}$ is bounded, this convergence still holds for $\{f(u_n)\}$ and $\{u_n^m\}$. Denote the limit function of $\{u_n\}$ as $u(t,x)$. Then $(u_n^m)_x=(v_n)_x\xrightarrow{w^\ast}(u^m)_x$ on any bounded domain of $[0,T]\times\mathbb R$ and it is easy to prove that $u(t,x)$ is continuous and is a weak solution to \eqref{0}. Obviously \eqref{zz0} and \eqref{pr4} imply that $u$ and $(u^m)_x$ are bounded, respectively.

\vskip 0.1in

\noindent (iv).At last, we will prove that $u(t,x)$ is a classic solution on $\{(t,x)\mid u(t,x)>0\}$.Suppose $u>0$ at a point $(t_0,x_0)$. Since $u$ is continuous, there exists a neighborhood $O\subset[0,T]\times\mathbb R$ and a constant $c>0$ such that $u(t,x)\geqslant c>0$ for any $(t,x)\in O$. Therefore, if $n$ is sufficiently large, it holds $u_n(t,x)\geqslant \frac{1}{2}c>0$ for any $(t,x)\in O$. Hence, $\{u_n\}$ is uniformly bounded and equicontinuous in $C^2(O)$. Thus, $u\in C^2(O)$ and satisfies the equation in classical sense.
\end{proof}

\vskip 0.2in

Next, we need $(u^{m-1})_x$ to be bounded. Let $u$ be a smooth positive classical solution of \eqref{0} in a rectangle $\Omega=(0,T_0]\times(a,b)$ and let $M_0=\max_{\overline{\Omega}}u$. Denote $\tilde v=u^{m-1}$, then $\tilde v$ satisfies
\begin{equation}\label{pro0}
\tilde v_t=-f^\prime(u)\tilde v_x+\frac{m}{m-1}\tilde v_x^2+m\tilde v\tilde v_{xx}
\end{equation}
in $\Omega$. We have the following Lemma.
\begin{lemma}\label{480}
Assume the condition in Lemma 4 holds. Let $\Omega^\ast=(\tau,T_0]\times(a_1,b_1)$ where $\tau>0$, $a_1>a$ and $b_1<b$, then
\begin{equation}\label{pro1}
|\tilde v_x(t,x)|\leqslant C(f,m,M_0,a_1-a,b-b_1,\tau)
\end{equation}
in $\overline{\Omega^\ast}$. If
\begin{equation*}
M_1\equiv\max_{[a,b]}\left|\frac{{\rm d}}{{\rm d}x}\big(u_0(x)^{m-1}\big)\right|<\infty,
\end{equation*}
then \eqref{pro1} holds in $[0,T_0]\times(a_1,b_1)$ and the constant $C$ now depends on $M_1$ instead of $\tau$.
\end{lemma}
\begin{proof}
Define
\begin{equation*}
G(r)=\frac{N}{3}r(4-r)
\end{equation*}
for $0\leqslant r\leqslant1$, where $N=M_0^{m-1}$. Then
\begin{equation}\label{pro2}
0\leqslant G\leqslant N,\quad \frac{2}{3}N\leqslant G^\prime\leqslant\frac{4}{3}N,\quad G^{\prime\prime}=-\frac{2}{3}N,\quad \left|\frac{G^{\prime\prime}}{G^\prime}\right|\leqslant1,\quad \left(\frac{G^{\prime\prime}}{G^\prime}\right)^\prime\leqslant-\frac{1}{4}.
\end{equation}
Since $0<\tilde v\leqslant N$, we can define a function $w(t,x)$ by $\tilde v=G(w)$. Then $0<w\leqslant1$ and the smoothness of $u$ carries over to $\tilde v$ and hence to $w$. It follows from \eqref{pro0} that in $\Omega$
\begin{equation}\label{pro3}
w_t=-f^\prime(u)w_x+mG\frac{G^{\prime\prime}}{G^\prime}w_x^2+\frac{m}{m-1}G^\prime w_x^2+mGw_{xx}.
\end{equation}

Setting $\beta=w_x$, differentiating \eqref{pro3} with respect to $x$ and multiplying the resultant equation by $\beta$, we obtain
\begin{equation}\label{pro4}
\begin{aligned}
\frac{1}{2}(\beta^2)_t-&mG\beta\beta_{xx}
=\left(\frac{m^2}{m-1}G^{\prime\prime}+mG\left(\frac{G^{\prime\prime}}{G^\prime}\right)^\prime\right)\beta^4\\
&\quad +\left(\frac{m(m+1)}{m-1}G^\prime+2mG\frac{G^{\prime\prime}}{G^\prime}\right)\beta^2\beta_x
-f^\prime(u)\beta\beta_x-\frac{u^{2-m}}{m-1}f^{\prime\prime}(u)G^\prime \beta^3
\end{aligned}
\end{equation}
in $\Omega$. Let $\zeta(t,x)$ be a $C^2(\overline{\Omega})$ function such that $\zeta=1$ in $\Omega^\ast$, $\zeta=0$ on the lower and lateral boundaries of $\Omega$, and
\begin{equation*}
0\leqslant\zeta\leqslant1,\quad 0\leqslant\zeta_t\leqslant\frac{2}{\tau},\quad |\zeta_x|\leqslant2\max\{a_1-a,b-b_1\}.
\end{equation*}
Set $z=\zeta^2\beta^2$, then $z\in C^2(\Omega)$. At a point $(t_0,x_0)\in\Omega$ where $z$ attains a maximum it holds
\begin{equation}\label{pro5}
0=\frac{1}{2}z_x=\zeta^2\beta\beta_x+\zeta\zeta_x\beta^2
\end{equation}
and
\begin{equation*}
z_t-mGz_{xx}\geqslant0.
\end{equation*}
The last inequality implies
\begin{equation}\label{pro6}
\zeta^2\left(\frac{1}{2}(\beta^2)_t-mG\beta\beta_{xx}\right)
\geqslant-3mG\zeta_x^2\beta^2+mG\zeta\zeta_{xx}\beta^2-\zeta\zeta_t\beta^2
\end{equation}
by using Cauchy's inequality. Applying \eqref{pro4} and \eqref{pro5} in \eqref{pro6}, we have
\begin{equation}\label{pro7}
\begin{aligned}
&-\left(\frac{m^2}{m-1}G^{\prime\prime}
+mG\left(\frac{G^{\prime\prime}}{G^\prime}\right)^\prime\right)\zeta^2\beta^4\leqslant \big(3mG\zeta_x^2-mG\zeta\zeta_{xx}+\zeta\zeta_t+f^\prime(u)\zeta\zeta_x\big)\beta^2\\
&\qquad\qquad\qquad\qquad -\left(\frac{m(m+1)}{m-1}G^\prime\zeta_x+2mG\frac{G^{\prime\prime}}{G^\prime}\zeta_x
+\frac{u^{2-m}}{m-1}f^{\prime\prime}(u)G^\prime\zeta\right)\zeta \beta^3
\end{aligned}
\end{equation}
which holds at $(t_0,x_0)$. From \eqref{pro2}, \eqref{pro7} implies
\begin{equation*}
2\zeta^2\beta^4\leqslant C_1\beta^2+2C_2\zeta|\beta|^3
\end{equation*}
at $(t_0,x_0)$. Thus, by using Cauchy's inequality, it follows that
\begin{equation*}
z(t,x)\leqslant z(t_0,x_0)=\big(\zeta^2\beta^2\big)(t_0,x_0)\leqslant C_1+C_2^2\equiv C_3.
\end{equation*}
Therefore
\begin{equation*}
\max_{\overline{\Omega^\ast}}|\tilde v_x|\leqslant\frac{4}{3}N\max_{\overline{\Omega^\ast}}|w_x|\leqslant\frac{4}{3}N\sqrt{C_3},
\end{equation*}
which completes the proof of the first assertion of this lemma. The proof of the second assertion is similar in which the main difference is to take $\zeta=\zeta(x)\in C_0^2([a,b])$ with $\zeta=1$ on $[a_1,b_1]$ and $0\leqslant\zeta\leqslant1$, so we omit the details.
\end{proof}

\vskip 0.2in

At last, we will prove the regularity of the solution.
\begin{lemma}\label{49}
Suppose $u_0(x)$ is continuous in $\mathbb R$, $0\leqslant u_0(x)\leqslant M$, $u_0^m,u_0^{m-1}$ are Lipschitz continuous, and $u(t,x)$ is a weak solution to the Cauchy problem \eqref{0}. Then\\
{\rm (1)}\quad $u\in C^{\frac{1}{2},1}([0,\infty)\times\mathbb R)$.\\
{\rm (2)}\quad $(u^m)_x$ exists and is continuous with respect to $x$. Especially, $(u^m)_x=0$ at the point where $u=0$.\\
{\rm (3)}\quad  $u_x$ exists and is continuous with respect to $x$. Especially, $u_x=0$ at the point where $u=0$.
\end{lemma}
The proof of Lemma \ref{49} is similar to the one of Theorem in \cite{aro69}, since the proof is based on Lemmas \ref{48} and \ref{480}, and independent of the equation itself, so we omit the details.

\vskip 0.2in

Thus the proof of Theorem \ref{thm1} is completed by Lemmas \ref{48} and \ref{49}.

\section{Time decay rate}
This section is devoted to the time-decay rate \eqref{jie1}. Firstly, we have
\begin{lemma}[Local estimate]\label{lemPhi}
Let $u(t,x)$ be the solution given in Theorem \ref{thm1} with $u(\tau,x)$ satisfying $\Phi(\tau,\cdot)\in H^1(\mathbb R)$ for any given $\tau\geqslant 0$, where $\Phi$ is defined in \eqref{phixx}. There exists $\Delta t>0$ independent of $\tau$ such that $\Phi(t,x)\in C\big(\tau,\tau+\Delta t;H^1(\mathbb R)\big)$ and
\begin{equation}\label{end1}
\sup_{\tau\leqslant t\leqslant\tau+\Delta t}\|\Phi(t,\cdot)\|_{H^1}\leqslant2\|\Phi(\tau,\cdot)\|_{H^1}.
\end{equation}
\end{lemma}
\begin{proof}
Multiplying \eqref{Phi}$_1$ by $\Phi$ and integrating the resultant equation, we have
\begin{equation}\label{91}
\frac{\textrm{d}}{\textrm{d}t}\int_\mathbb R|\Phi|^2\textrm{d}\xi+\int_\mathbb R|\Phi_\xi|^{m+1}\textrm{d}\xi\leqslant C\int_\mathbb R|\Phi|\Phi_\xi^2\textrm{d}\xi,
\end{equation}
where we have used Taylor's formula, \eqref{ab} and $f^{\prime\prime}>0,U^\prime\leqslant0$. Note that it holds $\Phi_\xi(t,\cdot)\in L^\infty(\mathbb R)$ for any $t>0$ from Lemma \ref{hhhh} and Theorem \ref{thm1}, we can conclude from \eqref{91} that
\begin{equation}\label{92}
\frac{\textrm{d}}{\textrm{d}t}\big(\|\Phi(t,\cdot)\|^2\big)\leqslant C\|\Phi(t,\cdot)\|^2
\end{equation}
by using
\begin{equation*}
\int_\mathbb R|\Phi|\Phi_\xi^2\textrm{d}\xi\leqslant\|\Phi_\xi\|_\infty^\frac{3-m}{2}\int_\mathbb R|\Phi||\Phi_\xi|^\frac{m+1}{2}\textrm{d}\xi\leqslant \frac{1}{2}\int_\mathbb R|\Phi_\xi|^{m+1}\textrm{d}\xi+C\int_\mathbb R\Phi^2\textrm{d}\xi.
\end{equation*}
Thus,
\begin{equation}\label{93}
\|\Phi(t,\cdot)\|\leqslant \|\Phi(\tau,\cdot)\|e^{C(t-\tau)}\leqslant2\|\Phi(\tau,\cdot)\|
\end{equation}
for any $t\in(\tau,\tau+\Delta t)$ by choosing $\Delta t$ suitably small.

\vskip 0.1in

We then need to estimate $\Phi_\xi$. Multiplying \eqref{Phi}$_1$ by $-\Phi_{\xi\xi}$, by a similar calculation, we have
\begin{equation*}
\frac{\textrm{d}}{\textrm{d}t}\left(\|\phi(t,\cdot)\|^2\right)\leqslant C\|\phi(t,\cdot)\|^2,
\end{equation*}
where we used the fact that $\phi_\xi(t,x)=\Phi_{\xi\xi}(t,x)\in L^\infty\big((0,\infty)\times\mathbb R\big)$ from Remark \ref{Urem} and Theorem \ref{thm1}. Thus,
\begin{equation}\label{94}
\|\phi(t,\cdot)\|\leqslant \|\phi(\tau,\cdot)\|e^{C(t-\tau)}\leqslant2\|\phi(\tau,\cdot)\|
\end{equation}
for any $t\in(\tau,\tau+\Delta t)$ by choosing $\Delta t$ suitably small. Comparing \eqref{93} and \eqref{94}, the proof is completed.
\end{proof}

\vskip 0.2in

Next, we will obtain the estimates of $\|\Phi(t,\cdot)\|_{H^1}$ on $t\in(0,T_1]$ for any $T_1$. That is,
\begin{lemma}[A priori estimate]\label{lemphi}
If $\|\Phi(t,\cdot)\|_{H^1}\leqslant2\varepsilon_0,t\in(0,T_1]$ for any $T_1>0$, it holds that
\begin{equation}\label{z1}
\begin{aligned}
&\|\Phi(t,\cdot)\|\leqslant\|\Phi_0\|,\\
&\|\phi(t,\cdot)\|\leqslant C\|\Phi_0\|_{H^1}^{\frac{1}{8}-\delta}(1+t)^{-\frac{1}{4(11m+7)}+\delta},
\end{aligned}
\end{equation}
where $\delta>0$ is any small constant and $C$ is independent of $t$ and $T_1$.
\end{lemma}
\begin{proof}
Multiplying \eqref{Phi} by $|\Phi|^{p-2}\Phi$, $p\geqslant2$, and using Taylor's expansion, we have
\begin{equation}\label{101}
\begin{aligned}
\frac{1}{p}\left(|\Phi|^p\right)_t+&(p-1)\big((U+\Phi_\xi)^m-U^m\big)|\Phi|^{p-2}\Phi_\xi
        -\frac{1}{p}f^{\prime\prime}(U)U^\prime|\Phi|^p\\
   &=(\cdots)_\xi-\frac{1}{2}f^{\prime\prime}(U+\theta_1\Phi_\xi)|\Phi|^{p-2}\Phi\Phi_\xi^2,
\end{aligned}
\end{equation}
where $\theta_1\in[0,1]$. Using \eqref{ab}, and noting that $f^{\prime\prime}>0,U^\prime\leqslant0$, it holds that
\begin{equation}\label{102}
\frac{\textrm{d}}{\textrm{d}t}\int_\mathbb R|\Phi|^p\textrm{d}\xi+\int_\mathbb R|\Phi|^{p-2}|\Phi_\xi|^{m+1}\textrm{d}\xi\leqslant C\int_\mathbb R|\Phi|^{p-1}\Phi_\xi^2\textrm{d}\xi.
\end{equation}
Thus, choosing $\varepsilon_0$ in \eqref{005} suitably small so that $\|\Phi\|$ is small, and hence, $\|\Phi\|_\lambda$ is small for any $2\leqslant\lambda<\infty$, we have
\begin{equation}\label{103}
\frac{\textrm{d}}{\textrm{d}t}\int_\mathbb R|\Phi|^p\textrm{d}\xi+\int_\mathbb R|\Phi|^{p-2}|\Phi_\xi|^{m+1}\textrm{d}\xi\leqslant0
\end{equation}
for $t\in(0,T_1]$ with some $T_1>0$ by using Proposition \ref{phii} and Lemma \ref{lemPhi}. Especially, we can let $T_1=\Delta t$ used in Lemma \ref{lemPhi}.
\begin{rem}
If we choose $p=2$ in \eqref{103}, it is easy to see that $\|\Phi_\xi\|_{m+1}^{m+1}\in L^1\big([0,T_1]\big)$ and $\|\Phi(t,\cdot)\|\leqslant\|\Phi_0\|$ for $t\in(0,T_1]$.
\end{rem}
Suppose $p>2$ and let $h(t)=\|\Phi(t,\cdot)\|_p^p$. Using Proposition \ref{phiii}, it holds
\begin{equation}\label{404a}
h^\prime+ch^\nu\leqslant0
\end{equation}
from \eqref{103}. Solving \eqref{404a} implies
\begin{equation}\label{405}
\|\Phi(t,\cdot)\|_p^p\leqslant C\|\Phi_0\|_p^p(1+t)^{-\frac{p-2}{3m+1}}
\end{equation}
for $t\in(0,T_1]$. Obviously, \eqref{405} also holds true for $p=2$.

\vskip 0.2in

Next we need the higher-order estimate. Multiplying \eqref{Phi} by $-|\Phi_\xi|^{q-2}\Phi_{\xi\xi}$ with $q\geqslant4$, and noticing that $|\Phi_\xi|^{q-2}\Phi_{\xi\xi}=\frac{1}{q-1}\left(|\phi|^{q-2}\phi\right)_\xi$, we have
\begin{equation}\label{201}
\begin{aligned}
\frac{1}{q(q-1)}\left(|\phi|^q\right)_t+&\frac{1}{q}f^{\prime\prime}(U)U^\prime|\phi|^{q}
    -\frac{1}{2}f^{\prime\prime}(U+\theta_1\phi)|\phi|^q\phi_\xi\\
&=(\cdots)_\xi-m\big((U+\phi)^{m-1}(U^\prime+\phi_\xi)-U^{m-1}U^\prime\big)|\phi|^{q-2}\phi_\xi.
\end{aligned}
\end{equation}
Noting $\phi\in L^\infty\big([0,\infty);L^\infty(\mathbb R)\big)$ by Proposition 3, integrating \eqref{201} with respect of $\xi$ over $\mathbb R$, we obtain
\begin{equation}\label{202}
\frac{\textrm{d}}{\textrm{d}t}\left(\|\phi\|_q^q\right)+mq(q-1)\int_\mathbb RB_1\textrm{d}\xi
\leqslant \mu\int_\mathbb R|\phi|^q\phi_\xi^2\textrm{d}\xi+C\int_\mathbb R|\phi|^q\textrm{d}\xi,
\end{equation}
where $\mu>0$ is a small constant, we have used Cauchy's inequality and the fact that $U,U+\phi\in L^\infty$ and $f\in C^2(\mathbb R)$, and
\begin{equation*}
B_1=\big((U+\phi)^{m-1}(U^\prime+\phi_\xi)-U^{m-1}U^\prime\big)|\phi|^{q-2}\phi_\xi.
\end{equation*}
Since $1<m<2$, the term $\int_\mathbb R|\phi|^q\phi_\xi^2\textrm{d}\xi$ can be majorized by some term like $\int_\mathbb R|\phi|^{m+q-3}\phi_\xi^2\textrm{d}\xi$ by choosing $\mu$ suitably small. In addition, $\int_\mathbb R|\phi|^q\textrm{d}\xi\leqslant C\int_\mathbb R|\phi|^{m+q-3}\textrm{d}\xi$. Then, we only need to deal with $\int_\mathbb RB_1\textrm{d}\xi$. In fact, we want to get the following inequality
\begin{equation}\label{203}
\frac{\textrm{d}}{\textrm{d}t}\left(\|\phi\|_q^q\right)+c\int_\mathbb R|\phi|^{m+q-3}\phi_\xi^2\textrm{d}\xi\leqslant C\int_\mathbb R|\phi|^{m+q-3}\textrm{d}\xi
\end{equation}
from \eqref{202}.

\vskip 0.2in

We will divide the integral $\int_\mathbb RB_1\textrm{d}\xi$ into several parts to discuss. Set
\begin{equation*}
\begin{aligned}
&D_0=[x_R,+\infty),\qquad D_1=\{\xi<x_R|\phi(\xi)\geqslant0\},\\
&D_2=\{\xi<x_R|\phi(\xi)<0,\phi_\xi(\xi)<0\},\\
&D_3=\{\xi<x_R|\phi(\xi)<0,\phi_\xi(\xi)\geqslant0\}.
\end{aligned}
\end{equation*}
Obviously, $\mathbb R=\cup_{i=0}^3D_i$ and $D_i\cap D_j=\emptyset(i\neq j)$ for any $i,j=0,1,2,3$.

\vskip 0.1in

\textbf{Part 1}. If $\xi\in D_0$, then $U=0$, $U^\prime=0$ a.e. and $\phi=u\geqslant0$. Thus, $B_1=\phi^{m+q-3}\phi_\xi^2$.

\vskip 0.1in

\textbf{Part 2}. For $\xi\in D_1$, noting $U>0$, we have
\begin{equation*}
\begin{aligned}
B_1&=(U+\phi)^{m-1}\phi^{q-2}\phi_\xi^2+\big((U+\phi)^{m-1}-U^{m-1}\big)U^\prime\phi^{q-2}\phi_\xi\\
   &\geqslant c(m)(U^{m-1}+\phi^{m-1})\phi^{q-2}\phi_\xi^2-C(m)|U^\prime|\phi^{m+q-3}|\phi_\xi|\\
   &\geqslant c(m)(U^{m-1}+\phi^{m-1})\phi^{q-2}\phi_\xi^2-\frac{c(m)}{2}\phi^{m+q-3}\phi_\xi^2-C(m)(U^\prime)^2\phi^{m+q-3}.
\end{aligned}
\end{equation*}
Then we have
\begin{equation*}
\int_{D_1}B_1\textrm{d}\xi\geqslant c\int_{D_1}\phi^{m+q-3}\phi_\xi^2\textrm{d}\xi-C\int_{D_1}\phi^{m+q-3}\textrm{d}\xi.
\end{equation*}

\vskip 0.1in

\textbf{Part 3}. If $\xi\in D_2$, choose a constant $0<C_1\ll1$ and define $C_2$ satisfying
\begin{equation*}
C_1=\frac{(1-C_2)^{m-1}}{1-(1-C_2)^{m-1}},
\end{equation*}
then $0\ll C_2(m)<1$. Let
\begin{equation*}
d=\frac{\phi}{U},\qquad B_2=\frac{(1+d)^{m-1}\phi_\xi+\big((1+d)^{m-1}-1\big)U^\prime}{(1+|d|^{m-1})\phi_\xi}.
\end{equation*}
It is easy to see that $|d|\leqslant1$ and $B_1=(U^{m-1}+|\phi|^{m-1})|\phi|^{q-2}\phi_\xi^2B_2$. We will then discuss case by case.

\vskip 0.1in

If $\left|\frac{U^\prime}{\phi_\xi}\right|>\frac{C_1}{2}$, then by Lemma 4,
\begin{equation*}
\begin{aligned}
B_1&=(U+\phi)^{m-1}|\phi|^{q-2}\phi_\xi^2+\big((U+\phi)^{m-1}-U^{m-1}\big)U^\prime|\phi|^{q-2}\phi_\xi\\
   &\geqslant|\phi|^{m+q-3}\phi_\xi^2-\frac{2}{C_1}|\phi|^{m+q-3}|U^\prime|^2-\frac{2}{C_1}\big|(U+\phi)^{m-1}
        -U^{m-1}\big||\phi|^{q-2}(U^\prime)^2\\
   &\geqslant|\phi|^{m+q-3}\phi_\xi^2-C(m,C_1,U^\prime)|\phi|^{m+q-3}.
\end{aligned}
\end{equation*}

\vskip 0.1in

If $|\frac{U^\prime}{\phi_\xi}|\leqslant\frac{C_1}{2}$ and $|d|\leqslant C_2$, we have $-1<d<0$, and
\begin{equation*}
\begin{aligned}
B_2&=\frac{1}{1+|d|^{m-1}}\left((1+d)^{m-1}\left(1+\left|\frac{U^\prime}{\phi_\xi}\right|\right)
        -\left|\frac{U^\prime}{\phi_\xi}\right|\right)\\
   &>\frac{1}{2}\left((1-C_2)^{m-1}\left(1+\left|\frac{U^\prime}{\phi_\xi}\right|\right)
        -\left|\frac{U^\prime}{\phi_\xi}\right|\right)\\
   &=\frac{1}{2}\left((1-C_2)^{m-1}-\big(1-(1-C_2)^{m-1}\big)\left|\frac{U^\prime}{\phi_\xi}\right|\right)\\
   &\geqslant\frac{1}{4}(1-C_2)^{m-1}
\end{aligned}
\end{equation*}
by $U^\prime<0$, so that
\begin{equation*}
B_1\geqslant\frac{1}{4}(1-C_2)^{m-1}\left(|\phi|^{m-1}+U^{m-1}\right)|\phi|^{q-2}\phi_\xi^2.
\end{equation*}

\vskip 0.1in

If $|\frac{U^\prime}{\phi_\xi}|\leqslant\frac{C_1}{2}$ and $C_2<|d|\leqslant1$, we have $|\phi|\thicksim U$ and $|\phi_\xi|\geqslant\frac{2}{C_1}|U^\prime|$. Since for any $\xi$ located in the left neighborhood of point $x_R$, $\phi<0$ and $\phi_\xi<0$ can not both be true, we can conclude, with the continuity of $U^\prime$ and $\phi_\xi$(from Lemma 4 in Section 3), and $\frac{2}{C_1}\gg1$, that this situation does not exist. In fact, if $|\phi_\xi|\geqslant\frac{2}{C_1}|U^\prime|$ at some point $\xi_1$, then by the continuity of $\phi_\xi$ and $U^\prime$, there exists a $\xi_2$ such that for any $\xi\in(\xi_1,\xi_2)$, $\phi_\xi\leqslant\frac{1}{C_1}U^\prime$. Integrating over $(\xi_1,\xi_2)$, it holds
\begin{equation*}
\phi(\xi_2)-\phi(\xi_1)\leqslant\frac{1}{C_1}\big(U(\xi_2)-U(\xi_1)\big).
\end{equation*}
Then, by using $|\phi|\thicksim U$ and $\phi<0$,
\begin{equation*}
U(\xi_2)+\phi(\xi_2)\leqslant\frac{1}{C_1}\big(U(\xi_2)-U(\xi_1)\big)+U(\xi_2)+\phi(\xi_1)<0,
\end{equation*}
which makes a contradiction with $U(\xi_2)+\phi(\xi_2)\geqslant0$.

\vskip 0.1in

\textbf{Part 4}. For $\xi\in D_3$, the discussion is similar to one in part 3.

\vskip 0.1in

If $|d|\leqslant C_2<1$, then $U\geqslant-\frac{1}{C_2}\phi$. By $\phi_\xi<0$ we get
\begin{equation*}
B_2=\frac{1}{1+|d|^{m-1}}\left((1+d)^{m-1}+\big((1+d)^{m-1}-1\big)\frac{U^\prime}{\phi_\xi}\right)
\geqslant\frac{1}{2}(1-C_2)^{m-1},
\end{equation*}
so it follows
\begin{equation*}
B_1\geqslant\frac{1}{2}(1-C_2)^{m-1}(U^{m-1}+|\phi|^{m-1})|\phi|^{q-2}\phi_\xi^2.
\end{equation*}

\vskip 0.1in

If $C_2<|d|\leqslant1$ and $\left|\frac{U^\prime}{\phi_\xi}\right|<C_1$, by a similar discussion in Part 3, this situation does not exist with the help of continuity of $\phi_\xi$ and $U^\prime$. In fact, if there exist $\xi_3<\xi_4$ such that for any $\xi\in(\xi_3,\xi_4)$, $\phi_\xi(\xi)>-\frac{1}{2C_1}U^\prime$, then integrating this inequality over $(\xi_3,\xi_4)$ implies
\begin{equation*}
\phi(\xi_4)-\phi(\xi_3)>\frac{1}{2C_1}\big(U(\xi_3)-U(\xi_4)\big).
\end{equation*}
Then, by using $|\phi|\thicksim U$ and $\phi<0$,
\begin{equation*}
U(\xi_3)+\phi(\xi_3)<U(\xi_3)+\phi(\xi_4)-\frac{1}{2C_1}\big(U(\xi_3)-U(\xi_4)\big)<0,
\end{equation*}
which makes a contradiction with $U(\xi_3)+\phi(\xi_3)\geqslant0$.

\vskip 0.1in

If $C_2<|d|\leqslant1$ and $\left|\frac{U^\prime}{\phi_\xi}\right|\geqslant C_1$, we have
\begin{equation*}
B_2>\frac{1}{2}\big(1-(1+d)^{m-1}\big)\frac{|U^\prime|}{\phi_\xi}>\frac{1}{4}C_1,
\end{equation*}
which means
\begin{equation*}
B_1>\frac{1}{4}C_1(U^{m-1}+|\phi|^{m-1})|\phi|^{q-2}\phi_\xi^2.
\end{equation*}

\vskip 0.1in

Now we can conclude that \eqref{203} holds true by the discussion from Part 1 to Part 4. Since $\phi_\xi\in L^\infty\big((0,\infty);L^\infty(\mathbb R)\big)$ by Lemma 4 in Section 3, we have, with the help of interpolation inequality, that
\begin{equation}\label{bu3}
\|\phi\|_\infty\leqslant C\|\phi_\xi\|_\infty^{\frac{p+1}{2p+1}}\|\Phi\|_p^{\frac{p}{2p+1}}\leqslant C\|\Phi\|_p^{\frac{p}{2p+1}}.
\end{equation}
Substituting \eqref{bu3} into \eqref{203} and using $\phi(t,\cdot)\in L^1(\mathbb R)$ and \eqref{405}, it holds
\begin{equation}\label{zz2}
\frac{\textrm{d}}{\textrm{d}t}\left(\|\phi\|_q^q\right)\leqslant C\|\phi\|_\infty^{m+q-4}\leqslant C\left(\|\Phi\|_p^p\right)^{\frac{m+q-4}{2p+1}}\leqslant C\|\Phi_0\|_p^{p\frac{m+q-4}{2p+1}}(1+t)^{-\frac{p-2}{3m+1}\frac{m+q-4}{2p+1}}
\end{equation}
for $t\in(0,T_1]$. Using Remark 2 and the H\"{o}lder continuity of $\phi$(see Lemma 4), we can prolong $\|\phi(t,\cdot)\|_q^q$ smoothly so that $\|\phi(t,\cdot)\|_q^q\in L^1\big([0,\infty)\big)$ and \eqref{zz2} holds on $(0,\infty)$. Choosing $p$ sufficiently large so that
\begin{equation}\label{zz3}
\frac{p-2}{3m+1}\frac{m+q-4}{2p+1}<2,
\end{equation}
and noting that $\phi\in L^\infty\big([0,\infty);L^\infty(\mathbb R)\big)$, then Lemma 3 implies
\begin{equation*}
\|\phi(t,\cdot)\|_q^q\leqslant C\|\Phi_0\|_p^{p\frac{m+q-4}{2p+1}}(1+t)^{-\frac{p-2}{2(3m+1)}\frac{m+q-4}{2p+1}}\leqslant C\|\Phi_0\|_{H^1}^{\frac{p-2}{2}\frac{m+q-4}{2p+1}}(1+t)^{-\frac{p-2}{2(3m+1)}\frac{m+q-4}{2p+1}}.
\end{equation*}
Using H\"{o}lder's inequality, we have
\begin{equation}\label{bu4}
\|\phi(t,\cdot)\|\leqslant \|\phi\|_1^{\frac{q-2}{2q-2}}\|\phi\|_q^{\frac{q}{2q-2}}\leqslant C\|\Phi_0\|_{H^1}^{\frac{p-2}{2}\frac{m+q-4}{2p+1}\frac{1}{2q-2}}(1+t)^{-\frac{p-2}{2(3m+1)}\frac{m+q-4}{2p+1}\frac{1}{2q-2}}.
\end{equation}
Let $p\rightarrow+\infty$ in \eqref{bu4}, and note that from \eqref{zz3}, $q<11m+8$, then Lemma \ref{lemphi} is proved.
\end{proof}

\vskip 0.2in

\begin{proof}[Proof of Theorem \ref{thm2}]
If we choose $\varepsilon_0$ suitably small and let $\tau=0$, we have, by comparing Lemmas \ref{lemPhi} and \ref{lemphi} and using \eqref{005}, that
\begin{equation}\label{z2}
\|\Phi(T_1,\cdot)\|_{H^1}\leqslant\frac{1}{2}\varepsilon.
\end{equation}
Then, let $\tau=T_1$, we can obtain \eqref{z2} with $T_1$ replaced by $2T_1$. Using similar analysis, \eqref{z1} holds for any $t>0$, so that \eqref{jie1} is proved.
\end{proof}
\vspace{0.3cm}
\textbf{Acknowledgment} \
The author would like to thank Prof. Feimin Huang for his helpful suggestions.

\end{document}